\documentclass[11pt]{amsart}
\usepackage{amssymb}
\usepackage{amsmath}
\usepackage{fancyhdr}
\usepackage[british]{babel}
\usepackage{geometry}
\usepackage{enumitem}
\usepackage{algpseudocode}
\usepackage{dsfont}
\usepackage{centernot}
\usepackage{xstring}
\usepackage{colortbl}
\usepackage[section]{algorithm}
\usepackage{graphicx}
\usepackage[font={footnotesize}]{caption}
\usepackage[usenames,dvipsnames,table]{xcolor}
\usepackage{tikz}
\usepackage[h]{esvect}
\usepackage[
		bookmarksopen=true,
		bookmarksopenlevel=1,
		colorlinks=true,
		linkcolor=darkblue,
        linktoc=page,
		citecolor=darkblue,
]{hyperref}
\usepackage{bbm}

\title[]{On a rainbow version of Dirac's theorem}
\date{\today}
\author[F.~Joos and J.~Kim]{Felix Joos and Jaehoon Kim}

\thanks{The research leading to these results was partially supported by the Deutsche Forschungsgemeinschaft (DFG, German Research Foundation) -- 339933727 (F.~Joos).
The research was also supported by the POSCO Science Fellowship of POSCO TJ Park Foundation and by the KAIX Challenge program of KAIST Advanced Institute for Science-X (J.~Kim).
}

\newtheorem{theorem}{Theorem}

\theoremstyle{definition}

\newtheorem{claim}{Claim}

\newtheoremstyle{stepstyle}{10pt}{5pt}{\em}{0pt}{\em}{:}{5pt}{}
\theoremstyle{stepstyle}

\definecolor{darkblue}{rgb}{0,0,0.5}

\def\noproof{{\unskip\nobreak\hfill\penalty50\hskip2em\hbox{}\nobreak\hfill%
       $\square$\parfillskip=0pt\finalhyphendemerits=0\par}\goodbreak}
\def\endproof{\noproof\bigskip}

\def\noclaimproof{{\unskip\nobreak\hfill\penalty50\hskip2em\hbox{}\nobreak\hfill%
       $-$\parfillskip=0pt\finalhyphendemerits=0\par}\goodbreak}

\def\endclaimproof{\noclaimproof\medskip}

\newdimen\margin
\def\textno#1&#2\par{
   \margin=\hsize
   \advance\margin by -4\parindent
          \setbox1=\hbox{\sl#1}
   \ifdim\wd1 < \margin
      $$\box1\eqno#2$$
   \else
      \bigbreak
      \hbox to \hsize{\indent$\vcenter{\advance\hsize by -3\parindent
      \it\noindent#1}\hfil#2$}
      \bigbreak
   \fi}

\def\lateproof#1{\removelastskip\penalty55\medskip\noindent\setcounter{claim}{0}\setcounter{step}{0}{\bf Proof of #1. }} % in each main proof, claim counter set back

\def\claimproof{\removelastskip\penalty55\medskip\noindent{\em Proof of claim: }}

\begin{document}

\newcommand{\new}[1]{\textcolor{red}{#1}}

\newcommand{\COMMENT}[1]{}
\renewcommand{\COMMENT}[1]{\footnote{\textcolor{blue!70!black}{#1}}} % comment out to hide comments
\newcommand{\TASK}[1]{}
\renewcommand{\TASK}[1]{\footnote{\textcolor{red!70!black}{#1}}} % comment out to hide comments

\newcommand{\todo}[1]{\begin{center}\textbf{to do:} #1 \end{center}}

\def\eps{{\varepsilon}}
\def\heps{{\hat{\varepsilon}}}

\newcommand{\ex}{\mathbb{E}}
\newcommand{\pr}{\mathbb{P}}
\newcommand{\cB}{\mathcal{B}}
\newcommand{\cA}{\mathcal{A}}
\newcommand{\cE}{\mathcal{E}}
\newcommand{\cS}{\mathcal{S}}
\newcommand{\cF}{\mathcal{F}}
\newcommand{\cG}{\mathcal{G}}
\newcommand{\bL}{\mathbb{L}}
\newcommand{\bF}{\mathbb{F}}
\newcommand{\bZ}{\mathbb{Z}}
\newcommand{\cH}{\mathcal{H}}
\newcommand{\cC}{\mathcal{C}}
\newcommand{\cM}{\mathcal{M}}
\newcommand{\bN}{\mathbb{N}}
\newcommand{\bR}{\mathbb{R}}
\def\O{\mathcal{O}}
\newcommand{\cP}{\mathcal{P}}
\newcommand{\cQ}{\mathcal{Q}}
\newcommand{\cR}{\mathcal{R}}
\newcommand{\cJ}{\mathcal{J}}
\newcommand{\cL}{\mathcal{L}}
\newcommand{\cK}{\mathcal{K}}
\newcommand{\cD}{\mathcal{D}}
\newcommand{\cI}{\mathcal{I}}
\newcommand{\cN}{\mathcal{N}}
\newcommand{\cV}{\mathcal{V}}
\newcommand{\cT}{\mathcal{T}}
\newcommand{\cU}{\mathcal{U}}
\newcommand{\cX}{\mathcal{X}}
\newcommand{\cZ}{\mathcal{Z}}
\newcommand{\cW}{\mathcal{W}}
\newcommand{\1}{{\bf 1}_{n\not\equiv \delta}}
\newcommand{\eul}{{\rm e}}
\newcommand{\Erd}{Erd\H{o}s}
\newcommand{\cupdot}{\mathbin{\mathaccent\cdot\cup}}
\newcommand{\whp}{whp }
\newcommand{\bX}{\mathcal{X}}
\newcommand{\bV}{\mathcal{V}}
\newcommand{\hbX}{\widehat{\mathcal{X}}}
\newcommand{\hbV}{\widehat{\mathcal{V}}}
\newcommand{\hA}{\widehat{A}}
\newcommand{\tA}{\widetilde{A}}
\newcommand{\hX}{\widehat{X}}
\newcommand{\hV}{\widehat{V}}
\newcommand{\tX}{\widetilde{X}}
\newcommand{\tV}{\widetilde{V}}
\newcommand{\cbI}{\overline{\mathcal{I}^\alpha}}
\newcommand{\hAj}{\widehat{A}^\sigma_j}
\newcommand{\hVM}{V^M}
\newcommand{\Haux}{H^{+}}
\newcommand{\Abad}{A^{bad}}
\newcommand{\Agood}{A'}
\newcommand{\Aupd}{A^{\sigma}}
\newcommand{\Agoodupd}{A^\ast}
\newcommand{\Gbad}{G^{bad}}
\newcommand{\Ggood}{G'}
\newcommand{\bfF}{\mathbf{F}}
\newcommand{\bfG}{\mathbf{G}}
\newcommand{\IND}{\mathbbm{1}}

\newcommand{\doublesquig}{%
  \mathrel{%
    \vcenter{\offinterlineskip
      \ialign{##\cr$\rightsquigarrow$\cr\noalign{\kern-1.5pt}$\rightsquigarrow$\cr}%
    }%
  }%
}

\newcommand{\defn}{\emph}

\newcommand\restrict[1]{\raisebox{-.5ex}{$|$}_{#1}}

\newcommand{\prob}[1]{\mathrm{\mathbb{P}}\left[#1\right]}
\newcommand{\expn}[1]{\mathrm{\mathbb{E}}\left[#1\right]}
\def\gnp{G_{n,p}}
\def\G{\mathcal{G}}
\def\lflr{\left\lfloor}
\def\rflr{\right\rfloor}
\def\lcl{\left\lceil}
\def\rcl{\right\rceil}

\newcommand{\qbinom}[2]{\binom{#1}{#2}_{\!q}}
\newcommand{\binomdim}[2]{\binom{#1}{#2}_{\!\dim}}

\newcommand{\grass}{\mathrm{Gr}}

\newcommand{\brackets}[1]{\left(#1\right)}
\def\sm{\setminus}
\newcommand{\Set}[1]{\{#1\}}
\newcommand{\set}[2]{\{#1\,:\;#2\}}
\newcommand{\krq}[2]{K^{(#1)}_{#2}}
\newcommand{\ind}[1]{$\mathbf{S}(#1)$}
\newcommand{\indcov}[1]{$(\#)_{#1}$}
\def\In{\subseteq}
\newcommand{\din}{d^{-}}

\begin{abstract} 
For a collection $\mathbf{G}=\{G_1,\dots, G_s\}$ of not necessarily distinct graphs on the same vertex set $V$, 
a graph $H$ with vertices in $V$ is a $\bfG$-transversal if there exists a bijection $\phi:E(H)\rightarrow [s]$ 
such that $e\in E(G_{\phi(e)})$ for all $e\in E(H)$.
We prove that for $|V|=s\geq 3$ and $\delta(G_i)\geq s/2$ for each $i\in [s]$, 
there exists a $\bfG$-transversal that is a Hamilton cycle. 
This confirms a conjecture of Aharoni.
We also prove an analogous result for perfect matchings.
\end{abstract}

\maketitle
\section{Introduction}

Suppose that we are given a collection $\mathbf{F}=\{F_1,\dots, F_s\}$ of not necessarily distinct subsets of some finite set $\Omega$.
Then a set $X\In \Omega$ such that $X\cap F_i\neq \emptyset$ for each $i\in [s]$ is often called a `transversal' of $\bfF$ or a `colourful' object of $\mathbf{F}$.
In the case where $\bfF$ is the edge set of a hypergraph, $X$ is known as a hypergraph transversal.
If $X=\{x_1,\ldots,x_s\}$ and $x_i\in F_i$ for all $i\in[s]$, then~$X$ is also called system of distinct representatives.
Frequently, we seek transversals with certain additional properties
as for example $|X\cap F_i|=1$ for all $i\in [s]$.

%An easy example to illustrate such a question is the following. Given an $n$-dimensional vector space $V$ and $n$ bases $B_1,\ldots,B_n$ of $V$, one can easily find a new basis $\{v_1,\dots, v_n\}$ of $V$ with $v_i \in B_i$ for each $i\in [n]$.

Other results that deal with transversals include results regarding transversals on Latin squares, 
a colourful version of Carath\'eodory's theorem by B\'ar\'any \cite{Barany:82}  and a further generalization by Holmsen, Pach and Tverberg~\cite{HPT:08}, 
a colourful version for a topological and a matroidal extension of Helly's theorem by Kalai and Meshulam~\cite{KM:05} and 
a colourful version of the Erd\H{o}s-Ko-Rado theorem by Aharoni and Howard~\cite{AH:17}.

Surprisingly, the study of `transversals' over collections of graphs has not received much attention until recently (for results on this topic see for example~\cite{ADGMS:18,Magnant:15}).
Here, we simply take~$\Omega$ to be the edge set of the complete graph on some vertex set $V$,
the set $\bfF$ as a collection of (the edge sets of) graphs with vertex set $V$,
and we ask for transversals (which are then collections of edges) with certain graph properties.
%Given connected graphs $G_1,\ldots,G_{n-1}$ on the same vertex set $V$ of size $n$,
%we can always find a spanning tree $T$ in $\bigcup_{i\in[n-1]}G_i$ and a bijection $\phi\colon E(T)\to [n-1]$ such that $e\in E(G_{\phi(e)})$ for each $e\in E(T)$.

%One may wonder whether this phenomena appears more often.
To be more precise, we define the following concept of transversals over a graph collection.
Let $\bfG=\{G_1,\dots, G_s\}$ be a collection of not necessarily distinct graphs with common vertex set~$V$. 
We say that a graph $H$ with vertices in $V$ is a \emph{partial $\bfG$-transversal} if there exists an injection $\phi \colon E(H)\to [s]$ such that $e\in E(G_{\phi(e)})$ for each $e\in E(H)$.
If in addition $|E(H)|=s$, then $H$ is a \emph{$\bfG$-transversal} (and $\phi$ a bijection).
We also say that $H$ is a path/cycle/triangle/matching (partial) $\bfG$-transversal if $H$ is a path/cycle/triangle/matching and similarly for other graphs.

% We call a collection $\cP$ of graphs \emph{condition}.  Then what condition $\cP$ of $n$-vertex graphs makes the following question true?

Let us consider the following question.
\begin{equation*}%\label{eq:metaquestion}
\begin{minipage}[c]{0.9\textwidth}\em
Let $H$ be a graph with $s$ edges, $\cG$ be family of graphs and $\bfG=\{G_1,\ldots, G_s\}$ be a collection of not necessarily distinct graphs on the same vertex set $V$ such that $G_i\in \cG$ for all $i\in [s]$.
Which properties imposed on $\cG$ yield a $\bfG$-transversal isomorphic to $H$?
\end{minipage}\ignorespacesafterend 
\end{equation*}

By considering the case when $G_1=\dots =G_s$, 
we need to study properties for $\cG$ such that~$H$ is a subgraph of each graph in $\cG$. 
However, this alone is not sufficient.
To see that, let $|V|=s\geq 5$ and $\cG$ be the collection of cycles with vertex set~$V$.
Consider $s-1$ identical cycles $G_1,\dots, G_{s-1}$ and another cycle $G_s$ which is edge-disjoint from the others. 
Then there does not exist a Hamiltonian $\bfG$-transversals; that is, one that is a Hamilton cycle (on $V$).

Similarly, consider the case when $G_1=\dots=G_{s-1}$ is a perfect matching in a cycle of length $2s$ and $G_s=\dots=G_{2s-2}$ is another perfect matching in the cycle. This shows that when $G_1,\dots, G_{2s-2}$
 are matchings of size $s$, a matching partial $\bfG$-transversal of size $s$ may not exist. It is conjectured that $2s-1$ matchings of size $s$ guarantees a matching partial $\bfG$-transversal of size $s$.
  Aharoni and Berger \cite{AB:09} improved a result of Drisko \cite{D:98} by verifying that if $\bigcup_{i\in [s]}G_i$ is a bipartite graph, then one can find a matching partial $\bfG$-transversal of size $s$ when $2s-1$ matchings of size $s$ are given. For non-bipartite case, Aharoni, Berger, Chudnovsky, Howard and Seymour \cite{ABCHS:19} showed that one can find a matching partial $\bfG$-transversal of size $s$  when $3s-2$ matchings of size $s$ are given. If we instead not only assume that $G_i$ is a matching but assume that $G_i$ satisfies the tight sufficient condition for a perfect matching, namely $\delta(G_i)\geq |G_i|/2=s/2$, then Theorem~\ref{thm: pm} shows that there is a $\bfG$-transversal that is a perfect matching when only $s$ graphs are given. Note that it is conjectured in \cite{AB:09} that $s$ matchings of size $s+2$ have a matching $\bfG$-transversal of size $s$.

Also it is not sufficient to impose the Tur\'an condition on the number of edges. 
In~\cite{ADGMS:18} (see also~\cite{Magnant:15}), it is shown that there is a triple of $n$-vertex graphs $G_1, G_2, G_3$ each having more than $n^2/4$ edges with no triangle transversal.
In fact, one needs to require (roughly) at least $0.2557n^2$  edges in each $G_i$ to guarantee the existence of a triangle transversal.

On the other hand,
Aharoni~\cite{ADGMS:18} conjectured that Dirac's theorem~\cite{Dir:52} can be extended to a colourful version and 
here we confirm this conjecture.

\begin{theorem}\label{thm: conj Aharoni}
Let $n\in \bN$ and $n\geq 3$.
Suppose $\bfG= \{G_1,\ldots,G_n\}$ is a collection of not necessarily distinct $n$-vertex graphs with the same vertex set
such that $\delta(G_i)\geq n/2$ for each $i\in[n]$.
Then there exists a Hamiltonian $\bfG$-transversal.
\end{theorem}
For the same reason as the bound in Dirac's theorem is sharp, 
we cannot improve upon the minimum degree bound in Theorem~\ref{thm: conj Aharoni}.
Cheng, Wang and Zhao~\cite{WY:19} recently proved a weaker version of Theorem~\ref{thm: conj Aharoni} with the condition $\delta(G_i)\geq (1/2+o(1))n$.

We also prove the following theorem concerning perfect matchings.

\begin{theorem}\label{thm: pm}
Let $n\in \bN$ and $n\geq 2$ even.
Suppose $\bfG= \{G_1,\ldots,G_{n/2}\}$ is a collection of not necessarily distinct $n$-vertex graphs with the same vertex set
such that $\delta(G_i)\geq n/2$ for each $i\in[n/2]$.
Then there exists a $\bfG$-transversal that is a perfect matching.
\end{theorem}

If all graphs $G_i$ above lie in a common balanced bipartite graph, then one can further improve the degree condition.
Aharoni, Georgakopoulos and Spr\"ussel \cite{AGS:09} proved a theorem regarding perfect matchings in $r$-partite $r$-uniform hypergraph. By considering an edge $uv\in G_i$ as a hyperedge $\{u,v,i\}$ and applying their theorem to the resulting $3$-uniform hypergraph, one can show that if $\bigcup_{i\in [n/2]} G_i$ is a balanced bipartite graph, then the condition $\delta(G_i)> n/4$ for each $i\in [n/2]$ is sufficient to find a $\bfG$-transversal that is a perfect matching.

\section{The proofs}

We write $[n]=\{1,\ldots,n\}$ and $[m,n] = \{m,m+1,\dots, n\}$. 
We denote by $\delta(G)$ the minimum degree of a graph~$G$.
For a digraph~$D$, we let~$A(D)$ be the arc set of $D$, and~$d^-_D(x)$ and~$d^+_D(x)$ refer to the indegree and outdegree of a vertex $x\in V(D)$, respectively.
We denote by~$N_D^-(x)$ the in-neighbourhood of $x\in V(D)$.

It will be also useful to specify a particular injection/bijection for a (partial) $\bfG$-transversal. 
To this end, we say that $(H,\phi)$ is a partial $\bfG$-transversal if $\phi\colon E(H)\to [s]$ 
satisfies $e\in E(G_{\phi(e)})$ for all $e\in E(H)$ and $\phi$ is injective
and a $\bfG$-transversal if $\phi$ is in addition also bijective.
If $i\notin \phi(E(H))$ for some $i\in[s]$, we say~$i$ is \defn{missed} by~$\phi$ and~$\phi$ \defn{misses}~$i$.

\lateproof{Theorem~\ref{thm: conj Aharoni}}
Assume for a contradiction that there do not exist Hamiltonian $\bfG$-transversals. 
It is routine to check the statement for $n\in \{3,4\}$, so we may assume that $n\geq 5$.
Let $V$ be the common vertex of the graphs in $\bfG$. 
For each $e\in \binom{V}{2}$, let 
$$c(e):=\{i\in [n]: e\in E(G_i)\}.$$

\begin{claim}\label{Cl2}
There exists a partial $\bfG$-transversal that is a cycle of length $n-1$.
\end{claim}
\claimproof 
Let $(C,\phi)$ be a partial $\bfG$-transversal which has the largest number of edges among all paths and cycles.
Among cycles and paths with the same number of edges, we prefer cycles.

Suppose $C= (x_1,\dots, x_{\ell+1})$ is an $\ell$-edge path with $\ell \in [3,n-1]$ (it is easy to see that $\ell\geq 3$ as $n\geq 5$ by simply picking the edges of $C$ greedily). 
Consider the $(\ell-1)$-edge path $P=(x_1,\dots, x_\ell)$.
The partial $\bfG$-transversal given by $\phi$ restricted to $E(P)$ misses at least two integers, say, $1$ and~$2$. 
Then $1,2\notin c(x_1x_\ell)$, as otherwise $(x_1,\dots, x_\ell,x_1)$ forms an $\ell$-edge cycle partial $\bfG$-transversal which contradicts the choice of~$(C,\phi)$.
Let 
\begin{align*}
 I_1:=\{i\in [\ell-2]\colon 1\in c(x_1x_{i+1})\}\text{ and } 
 I_2:=\{i\in [2,\ell-1] \colon 2\in c(x_i x_\ell)\}.
 \end{align*}
 Note that we have
\begin{align}\label{eq: no y}
|N_{G_1}(x_1) \setminus V(P)| + |N_{G_2}(x_\ell) \setminus V(P)| 
\leq  n - \ell,
\end{align}
 otherwise, by the pigeonhole principle, there exists $y\in V\setminus V(P)$ such that $(x_1,\dots, x_\ell,y,x_1)$ forms an $(\ell+1)$-edge cycle partial $\bfG$-transversal, 
 again a contradiction to the choice of~$(C,\phi)$. 
 Since $\delta(G_i)\geq n/2$ for all $i\in [n]$ and $1,2\notin c(x_1x_\ell)$, equation \eqref{eq: no y} implies that
\begin{align*}
|I_1| + |I_2| 
\geq n/2+ n/2 - |N_{G_1}(x_1) \setminus V(P)| - |N_{G_2}(x_\ell) \setminus V(P)| 
\geq \ell.
\end{align*}
As $I_1\cup I_2\subseteq [\ell-1]$, there exists an integer $j \in I_1\cap I_2 \subseteq [2,\ell-2]$.
Hence deleting ${x_jx_{j+1}}$ from~$E(P)$ and adding $x_1x_{j+1},x_jx_\ell$ yields a partial $\bfG$-transversal that is a cycle of length~$\ell$,
which is a contradiction to the choice of $(C,\phi)$. 
Hence we may assume that~$C$ is a cycle.

In view of the statement, we may assume that $C=(x_1,\dots, x_\ell,x_1)$ is an $\ell$-edge cycle for some $\ell\in [3,n-2]$ and there are two integers, say $1$ and $2$, that are missed by $\phi$.
Observe that
$\ell\geq n/2+1$, since otherwise we have
$$ |N_{G_1}(x_1) \setminus V(C)| \geq 1 \text{ and } |N_{G_2}(x_\ell) \setminus V(C)|  \geq 1,$$
and we obtain two not necessarily distinct vertices $y,z \in V\setminus V(C)$ with $1\in c(x_1y), 2\in c(x_\ell z)$. 
Then $(y, x_1,\dots, x_\ell, z)$ is a partial $\bfG$-transversal which is either path or cycle with $\ell+1$ edges and this contradicts the choice of $(C,\phi)$.

We claim that, for each $v\in V\setminus V(C)$ and $i\in [2]$, we have $N_{G_i}(v)\subseteq V(C)$.
Suppose not. Then there exists $i\in [2]$ and $u,v \in V\setminus V(C)$  with $uv \in E(G_i)$.
As we have $d_{G_{3-i}}(v) \geq n/2 > |V\setminus V(C)|$, we have, by symmetry, $x_\ell v \in E(G_{3-i})$.
Consequently, $(x_1,\dots, x_\ell, v, u)$ contradicts the choice of~$(C,\phi)$.
Thus, for each $v\in V\setminus V(C)$ and $i\in[2]$, we have $N_{G_i}(v)\subseteq V(C)$.

Fix some $v\in V\setminus V(C)$.
Let
$$I_1:=\{i\in [\ell]\colon 1\in c(vx_{i+1})\} \text{ and } I_2:=\{i\in [\ell] \colon 2\in c(v x_i)\},$$ 
where we identify $x_{\ell+1}$ with $x_1$.
Then 
$$|I_1| + |I_2| \geq \delta(G_1) + \delta(G_2) \geq n > \ell,$$
and there exists an integer $j\in I_1\cap I_2$.
Hence deleting ${x_jx_{j+1}}$ from $E(C)$ and adding $vx_{j},vx_{j+1}$ yields partial $\bfG$-transversal that is a cycle of length $\ell+1$, which is a contradiction to the choice of $(C,\phi)$.
This proves Claim~\ref{Cl2}.
\endclaimproof

\medskip

By Claim~\ref{Cl2}, there exists a cycle partial $\bfG$-transversal $(C,\phi)$ with $C:=(x_1,\dots,x_{n-1},x_1)$.
By relabelling colours, we may assume that $\phi(x_ix_{i+1})=i$ for each $i\in [n-1]$ where we identify~$x_n$ with~$x_1$.
Hence~$\phi$ misses $n$.
Let $\{y\}=V\setminus V(C)$.
We consider the following auxiliary digraph~$D$ on vertex set $[n]$ such that
$$A(D) = \bigcup_{i\in [n-1]} \{ x_i z : z\neq x_{i+1}, i\in c(x_i z)  \}.$$
As $\delta(G_i)\geq n/2$ for all $i\in [n-1]$ and thus $d^+_D(x)\geq n/2-1$ for all $x\in V(C)$, we obtain that 
$|A(D)|  \geq (n-1)({n}/{2}-1)$.
Let $$I:=\{i\in[n-1]\colon x_iy \in A(D)\} \text{ and } I':= \{i\in[n-1]\colon x_{i+1}y \in E(G_n)\}.$$
We claim that $\din_{D}(y) \leq \frac{n}{2}-1$.
Otherwise, 
we have $|I| + |I'| \geq \din_{D}(y) + \delta(G_n) > n-1 = |V(C)|$.
So, there exists $j\in I\cap I'$ and thus $(E(C)\setminus \{x_jx_{j+1}\})\cup \{x_jy,yx_{j+1}\}$ is the edge set of a Hamiltonian $\bfG$-transversal,
which is a contradiction.
 
Hence, we assume from now on that $\din_{D}(y) \leq \frac{n}{2}-1$. 
By our definition of $D$, we have $d^+_{D}(y) = 0$ and thus
\begin{align}\label{eq: AD}
|A(D-y)|
\geq (n-1)\left(\frac{n}{2}-1\right) - \frac{n}{2}+1
%=(n-1)\left(\frac{n}{2}-1- \frac{1}{2}\cdot \frac{n-2}{n-1}\right)
> 
(n-1)\left(\frac{n}{2}-\frac{3}{2}\right).
\end{align}
%Now, we consider the following two cases.\COMMENT{Old proof didn't consider the possibility that $x^{k+1}=y$. In such a case, old proof dose not work. So we need to divide these two cases.
%} \newline

%\noindent {\bf Case 1.} There is a vertex, say $x_1$, such that $\din_{D-y}(x_1) > n/2-1$.

Let us assume for now that there exists a vertex, say $x_1$, such that $\din_{D}(x_1) > n/2-1$.
%Let $\hat{I}:=\{i\in [2, n-2]\colon i \in c(x_1x_i)\}$. 
Consequently, we conclude that
\begin{align}\label{eq:y1}
|\{i\in [2, n-2]\colon i \in c(x_1x_i)\}| = \din_{D}(x_1) \geq \frac{n}{2}- \frac{1}{2}.
\end{align}
Let 
$$I_1:= \{i\in[n-1]\colon x_{i}y \in E(G_1)\} \text{ and } I_n:= \{i\in[n-1]\colon x_{i+1}y \in E(G_n)\}.$$
Clearly, $|I_1|+|I_n|\geq n$, so there exists a $j\in I_1\cap I_n$.
We may assume that $j\neq 1$ as otherwise $(E(C)\setminus \{x_1x_2\})\cup \{x_1y,x_2y\}$ is the edge set of a Hamiltonian $\bfG$-transversal,
which is a contradiction.

Let $(P,\phi')$ with $P= (x_2,\ldots, x_j, y, x_{j+1},\ldots ,x_{n-1},x_1)$ be a path partial $\bfG$-transversal that arises from~$\phi$ by deleting $\{x_1x_2,x_jx_{j+1}\}$ from its domain and by setting $\phi'(x_jy):=1$ and $\phi'(x_{j+1}y):=n$. 
Observe that $\phi'$ misses (only) $j$.
We write $P= (x^1,\ldots, x^n)$ such that $x^1=x_2$.
Let 
$$J_1:=\{i\in [n-2]\colon j\in c(x^1x^{i+1})\} \text{ and } J_n:=\{ i\in [n-2] \colon x^i \in N^{-}_{D}(x_1) \}.$$
If $j\in c(x^1x^{n})$, then there is a Hamiltonian $\bfG$-transversal, which is a contradiction; so $|J_1| \geq \delta(G_j)$.
Also, as $x^{n-1}\in \{x_{n-1},y\}$,
the definition of $D$ ensures that $x^{n-1}\notin N^{-}_{D}(x_1)$. 
Hence \eqref{eq:y1} implies that $|J_n|\geq  n/2 -1/2$ and thus $|J_1|+|J_2| \geq n$.
Since $J_1\cup J_2 \subseteq [n-2]$, 
there exist at least two integers in $J_1\cap J_n$
and at least one of them, say $k$, satisfies $x^{k+1} \neq y$.
Moreover, $x^k \neq y$ as $y\notin N^{-}_{D}(x_1)$. 
Hence, $\phi'(x^k x^{k+1})\in c(x^kx^n)$ and  
$(E(P)\setminus \{x^kx^{k+1}\}) \cup \{x^1x^{k+1},x^kx^n\}$ forms a Hamiltonian $\bfG$-transversal, which is a contradiction. 

\medskip

Therefore, we may assume that $\din_{D}(x_i) \leq {n}/{2}-1$ for all $i\in [n-1]$.
We define 
$$\cJ:= \left\{ i\in [n-1]: \din_{D}(x_i) = \left\lfloor \frac{n}{2}-1 \right\rfloor \right\}.$$
Then \eqref{eq: AD} implies that 
$$ \left\lfloor \frac{n}{2}-1 \right\rfloor |\cJ| + \left\lfloor \frac{n}{2}-2 \right\rfloor( n-1-|\cJ|) 
\geq |A(D)| > (n-1) \left(\frac{n}{2}-\frac{3}{2}\right).$$
Hence, we have
$$|\cJ|  
>(n-1)\left( \frac{n}{2} - \left\lfloor \frac{n}{2} \right\rfloor +\frac{1}{2} \right)
\geq  \frac{n-1}{2}.$$
Let $\cJ':= \{i\in[n-1]\colon x_{i+1}y \in E(G_n)\}$.
Clearly, $|\cJ|+|\cJ'|\geq n$ and so there exists a $j\in \cJ \cap \cJ'$.
Let $(Q,\phi')$ with $Q=(y, x_{j+1}, x_{j+2},\dots, x_{n-1}, x_1, \dots, x_j)$ be a path partial $\bfG$-transversal that arises from~$\phi$ by deleting $\{x_jx_{j+1}\}$ from its domain and by setting $\phi'(x_{j+1}y):=n$. 
Observe that~$\phi'$ misses $j$.
We write $Q=(x^1,\ldots,x^n)$ such that $x^1=y$.
Let 
$$J_1:=\{i\in [n-2]\colon j\in c(x^1x^{i+1})\} \text{ and } J_n:=\{ i\in [2,n-2] \colon x^i \in N^{-}_{D}(x^n) \}.$$
If $j\in c(x^1x^n)$, then there is a Hamiltonian $\bfG$-transversal, which is a contradiction; so $|J_1| \geq \delta(G_j) \geq n/2.$
Note that $x^1 =y \notin N^{-}_{D}(x^n)$ and $x^{n-1}= x_{j-1} \notin N^{-}_{D}(x^n)$ by the definition of~$D$.
As $x^n =x_j \in \cJ$, we infer that
$|J_n| = \lfloor n/2-1\rfloor$.
We obtain $|J_1|+|J_n| \geq n-1$.
As $J_1\cup J_n \subseteq [n-2]$, 
there exists an integer $k\in J_1\cap J_n \subseteq [2,n-2]$.
Since $x^k \neq y=x^1$,
we conclude that
$\phi'(x^k x^{k+1})\in c(x^kx^n)$ and  
$(E(P)\setminus \{x^kx^{k+1}\}) \cup \{x^1x^{k+1},x^kx^n\}$ contains a  Hamiltonian $\bfG$-transversal.
This is the final contradiction.
\endproof

\lateproof{Theorem~\ref{thm: pm}}
We use similar notation as in the proof of Theorem~\ref{thm: conj Aharoni};
in particular, let $V$ be the common vertex set of the graphs in $\bfG$ and for each $e\in \binom{V}{2}$, let 
$$c(e):=\{i\in [n/2]: e\in E(G_i)\}.$$
For a partial $\bfG$-transversal $(M,\phi)$,
we refer to $|E(M)|$ as the \defn{size} of $(M,\phi)$.
We assume for a contradiction that there does not exist a matching $\bfG$-transversal.

It is easy to see that $G$ contains a matching partial $\bfG$-transversal of size $n/2-1$.
Indeed, consider a matching partial $\bfG$-transversal $(M,\phi)$ of maximum size $\ell$.
Assume for a contradiction that $\ell<n/2-1$ and $\phi$ misses 1 and 2, say.
Clearly, $\{1,2\}\cap c(xx')=\emptyset$ for all $xx'\in \binom{V\setminus V(M)}{2}$.
Fix two vertices $x,x'\in  V\setminus V(M)$.
Let the weight of an edge $e=uv$ be $\IND_{1\in c(xu)}+\IND_{1\in c(xv)}+\IND_{2\in c(x'u)}+\IND_{2\in c(x'v)}$.
Since $\delta(G_i)\geq n/2$ for $i\in [2]$, we deduce that the sum of the weights of the edges in $M$ is at least $n$.
Hence there is an edge $e=yy'$ in $M$ with weight at least~$3$.
Replacing $e$ by $\{xy,x'y'\}$ or $\{x'y,xy'\}$ yields a contradiction to our assumption that the size of $(M,\phi)$ is maximum.

For a contradiction, we assume that there is no matching $\bfG$-transversal. 
Let $\ell:=n/2-1$.
For a matching partial $\bfG$-transversal $(N,\phi)$, 
we let  $D_N^\phi$ be a digraph with vertex set~$V$ and 
\begin{align*}
A(D_N^\phi):=\{xy:\phi(xz)\in c(xy),y\neq z,xz\in E(N)\}.
\end{align*}

\begin{claim}\label{claim:indegree1}
	$d^-_{D^\phi_M}(x)\leq \ell-1$ for all matching partial $\bfG$-transversals $(M,\phi)$ of size $\ell$ and $x\in V\setminus V(M)$.
\end{claim}
\claimproof
We define $D:=D_M^\phi$.
We assume for a contradiction that $d^-_{D}(x)\geq \ell$.
Let $\{x'\}= V\setminus (V(M)\cup \{x\})$.
Say $\phi$ misses $1$.
Clearly, $1\notin c(xx')$.
As $\delta(G_1)\geq n/2$ and $d^-_{D}(x)\geq \ell$, there exists an edge $yy'\in V(M)$ such that $yx\in A(D)$ and $1\in c(x'y')$.
However, then removing~$yy'$ from $M$ and adding $xy$ and $x'y'$ yields a matching $\bfG$-transversal, which is a contradiction.
\endclaimproof

\begin{claim}\label{claim:indegree2}
	$d^-_{D_M^\phi}(x)\leq \ell$ for all matching partial $\bfG$-transversals $(M,\phi)$ of size $\ell$ and $x\in V$.
\end{claim}
\claimproof
We define $D:=D_M^\phi$.
We assume for a contradiction that $d^-_{D}(x)\geq \ell+1$ and~$\phi$ misses $1$, say.
By Claim~\ref{claim:indegree1}, we conclude that $x\in V(M)$.
Let $y$ be the neighbour of $x$ in $M$ and $\phi(xy)=2$, say.
Let $\{z,z'\}=V\setminus V(M)$.
Suppose $i\in c(y\tilde{z})$ for some $i\in [2],\tilde{z}\in \{z,z'\}$.
Then let $(M',\phi')$ be the matching partial $\bfG$-transversal where $(M',\phi')$ arises $(M,\phi)$ by deleting $xy$ from~$M$, adding $y\tilde{z}$, and
assigning~$i$ on $y\tilde{z}$.
Hence, for $D':=D_{M'}^{\phi'}$, we obtain $d^-_{D'}(x)\geq \ell+1$, which is a contradiction to Claim~\ref{claim:indegree1}.
So we may assume that $\{1,2\}\cap (c(yz)\cup c(yz'))=\emptyset$.

Let $V':= V\setminus \{x,y,z\}$.
Then $|N_{G_2}(y)\cap V'|\geq n/2 -1$ and $|N_M(N_{G_1}(z))\cap V'|\geq n/2 -1$.
Consequently, there exists a vertex $u\in V'\cap N_{G_2}(y) \cap N_M(N_{G_1}(z))$.
Observe that $u\notin \{x,y,z,z'\}$.
Let $u'$ be the neighbour of $u$ in $M$.
Let $(M'',\phi'')$ be the matching partial $\bfG$-transversal where~$M''$ arises~$M$ by deleting $xy,uu'$ and adding $uy,u'z$ 
and~$\phi''$ arises from~$\phi$ by assigning~$u'z$ to~$1$ and~$uy$ to~$2$.
We write $D''$ for $D_{M''}^{\phi''}$ and observe that $d^-_{D''}(x)\geq d^-_{D}(x)-1$
as $y\in N^-_{D''}(x)\setminus N^-_{D}(x)$ and $N^-_{D}(x)\setminus N^-_{D''}(x)\In \{u,u'\}$.
However, exploiting Claim~\ref{claim:indegree1}, $(M'',\phi'')$ yields a contradiction.
\endclaimproof

\begin{claim}\label{claim:indegree3}
	For all matching partial $\bfG$-transversals $(M,\phi)$ of size $\ell$,
	there are at least $n/2$ vertices $x\in V(M)$ with $d^-_{D_M^\phi}(x)\geq \ell-1$.
\end{claim}
\claimproof
We define again $D:=D_M^\phi$.
Observe that the number of arcs in $D$ is at least $2\ell^2$,
as $d^+_D(x)\geq \ell$ for all $x\in V(M)$.
Assuming that there are at most $n/2-1=\ell$ vertices $x\in V(M)$ with $d^-_D(x)\geq \ell-1$,
implies in view of Claims~\ref{claim:indegree1} and \ref{claim:indegree2} 
that $|A(D)|\leq \ell^2 + \ell(\ell-2) +2(\ell-1)<2\ell^2$, which is a contradiction.
\endclaimproof

Let $(M,\phi)$ be some matching partial $\bfG$-transversal of maximum size.
In view of the above, the size of~$M$ equals~$\ell$ and so~$\phi$ misses~$1$, say.
Let $\{z,z'\}=V\setminus V(M)$ and $D:=D_M^\phi$.
By Claim~\ref{claim:indegree3} and as $\delta(G_1)\geq n/2$,
there exists $xy\in V(M)$ with $d^-_D(x)\geq \ell-1$ and $1\in c(yz)$.
Say, $\phi(xy)=2$.
Let $(M',\phi')$ arise from $(M,\phi)$ by deleting $xy$ from $M$, adding $yz$ and assigning~$yz$ to~1.
Let $D':=D_{M'}^{\phi'}$. %and $V':=V\setminus \{x,y,z'\}$.

\begin{claim}\label{claim:non colours} The following hold:
	\begin{enumerate}[label=\rm  (\alph*)]
		\item\label{item:out x} $|N^-_{D'}(x)\cap (V\setminus\{x,z,z'\})|\geq \ell-1$;
		\item\label{item:z'} $|N_{G_2}(z')\cap (V\setminus\{x,y,z'\})|\geq n/2$.
	\end{enumerate}
\end{claim}
\claimproof
Statement~\ref{item:out x} is obvious.
To see~\ref{item:z'}, we first observe that if $2\in c(xz')$, then we can delete~$xy$ from~$M$ and add~$xz'$ and~$yz$ and obtain a matching $\bfG$-transversal.
Moreover, if $2\in c(yz')$, then the matching that arises from $M$ by deleting $xy$, adding $yz'$, and assigning $2$ to $yz'$ contradicts Claim~\ref{claim:indegree1}.
This proves~\ref{item:z'}.
\endclaimproof

Observe that $N_{G_2}(z')\cap (V\setminus\{x,y,z'\})\In V(M')$.
Let $A$ be the set of vertices that are joined by an edge in $M'$ to a vertex in $N_{G_2}(z')\cap (V\setminus\{x,y,z'\})$.
Consequently, $A\In V\setminus\{x,z,z'\}$ and $|A|\geq n/2$ by Claim~\ref{claim:non colours}\ref{item:z'}.
As $|V\setminus\{x,z,z'\}|=n-3< n/2 + \ell-1 \leq |A|+ |N^-_{D'}(x)\cap (V\setminus\{x,z,z'\})|$,
there is a vertex $u\in A\cap N^-_{D'}(x)\cap (V\setminus\{x,z,z'\})$.
Let $v$ be the neighbour of $u$ in $M'$.
Deleting $uv$ and adding $ux$ and $vz'$ to $M'$ gives rise to a matching $\bfG$-transversal.
This is the final contradiction and completes the proof.
\endproof

\bibliographystyle{amsplain}
\bibliography{rainbow_HC_lit}

\vfill
%
%\small
%\vskip2mm plus 1fill
%\noindent
%Version \today{}
%\bigbreak
%
%
%

%\noindent
%Felix Joos
%{\tt <f.joos@bham.ac.uk>}\\
%School of Mathematics\\ 
%University of Birmingham\\
%United Kingdom\\

%\noindent
%Jaehoon Kim
%{\tt <Jaehoon.Kim.1@warwick.ac.uk>}\\
%Mathematics Institute\\ 
%University of Warwick\\
%United Kingdom\\

\begin{minipage}{0.55\textwidth}
Felix Joos\\
{\tt <joos@informatik.uni-heidelberg.de>}\\
Institute for Computer Science\\ 
Universit{\"a}t Heidelberg\\
Germany\\
\end{minipage}%
%\hfill
\begin{minipage}{0.55\textwidth}
%\begin{tabular}{p{\textwidth}}
Jaehoon Kim\\
{\tt <jaehoon.kim@kaist.ac.kr>}\\
Department of Mathematical Sciences\\ 
KAIST\\
Republic of Korea\\
%\end{tabular}
\end{minipage}

\end{document}